\documentclass[10pt]{amsproc}

\newcommand\query[1]{}

\def\dis{\displaystyle}
\newcommand{\vp}{\varphi}
\newcommand{\cald}{{\mathcal D}}
\newcommand{\bbbn }{{\mathbb N}}

\DeclareMathOperator{\artA}{\Lambda}
\DeclareMathOperator{\ch}{char}
\DeclareMathOperator{\Der}{Der}

\DeclareMathOperator{\Grass}{Grass}
\DeclareMathOperator{\Hoch}{HH}
\DeclareMathOperator{\Inn}{Inn}
\DeclareMathOperator{\ord}{ord}

\theoremstyle{definition}
\newtheorem{definition}{Definition}[section]
\newtheorem{sit}[definition]{}
\newtheorem{remark}[definition]{Remark}

\newtheorem{example}[definition]{Example}

\theoremstyle{plain}
\newtheorem{proposition}[definition]{Proposition}
\newtheorem{theorem}[definition]{Theorem}
\newtheorem{lemma}[definition]{Lemma}
\newtheorem{corollary}[definition]{Corollary}

\begin{document}
\title[Loops but no Outer Derivations]
{Artin Algebras with Loops\\
but no Outer Derivations}

\author[R.-O.~Buchweitz]{Ragnar-Olaf~Buchweitz}
\address{Department of Mathematics, University of Toronto,
Toronto, Ontario, Canada M5S 3G3}
\email{ragnar@math.utoronto.ca}

\author[S.~Liu]{Shiping Liu}
\address{D\'epartement de math\'ematiques et d'informatique,
Universit\'e de Sherbrooke,
Sherbrooke, Qu\'ebec, Canada J1K 2R1}
\email{shiping.liu@dmi.usherb.ca}

\thanks{
Each author is partly supported by a grant from NSERC;
\endgraf
This collaboration was partly supported by an Ontario-Quebec Exchange
Grant
of the Ministry of Education and Training of Ontario, Canada.
}
\date{\today}
\maketitle

\section{The Result}

If $\artA $ is a finite dimensional algebra over an algebraically
closed field $k$, then there exists a unique finite quiver
$Q_{\artA }$, the {\it ordinary quiver} of $\artA$, such that
$\artA $ is Morita equivalent to a
quotient of the path algebra $kQ_{\artA }$ by an admissible ideal $I$; 
see, for example, \cite{Ga} or \cite{ARS}.

The first Hochschild cohomology group of $\artA $ over $k$ is
$$
\Hoch^{1}(\artA) = \Hoch^{1}(\artA/k) = \frac{\Der_{k}(\artA ,\artA
)}{\Inn(\artA )}\,,
$$
the $k$-vectorspace of {\em outer $k$-derivations\/}, that is, the
quotient of all $k$-derivations of $\artA $ modulo the inner ones.

It has been suspected for some time, and \cite{Ci} seems to be the 
earliest reference, that
vanishing of the first Hochschild cohomology precludes the existence 
of oriented cycles in the ordinary quiver.  There are various 
supporting partial results, among them \cite[(2.3), (3.2)]{Hap}, 
\cite[(2.2)]{BM}, \cite[(1.3)]{Gas}.

An algebra $\artA $ without oriented cycles in its ordinary quiver is 
of finite global dimension.  In turn, finite global dimension implies 
that there are no loops in the ordinary quiver, see \cite{Len} 
or \cite{Ig}.

The following result shows that even the existence of loops in the
ordinary quiver
is no guarantee for the existence of non trivial outer derivations,
thus refuting the above suspicion.

\begin{theorem}
\label{main result}
		Let $Q$ be the quiver with three vertices ${\bf 1,2,3}$ and 
		three arrows $\alpha\colon{\bf 2}\to {\bf 1},\ \beta\colon{\bf 
		2}\to {\bf 2},\ \gamma\colon{\bf 3}\to {\bf 2}$.  For any 
		field $k$ there exist finite dimensional $k$-algebras with 
		ordinary quiver $Q$ that admit only inner derivations.
		These algebras are necessarily of infinite global dimension 
		and of infinite representation type.
\end{theorem}

The proof will occupy the rest of this paper.  We first describe the 
basic structure of the finite dimensional algebras $\artA$ 
with ordinary quiver $Q$, then identify $\Hoch^{1}(\artA)$ as a space 
of first order differential operators in one variable, give 
in (\ref{exponential examples}) a quick 
proof for fields whose prime field is large enough
and present finally in (\ref{positive characteristic}) 
a family of examples that works over any field.

\section{The Basic Structure}

We begin by exhibiting the path algebra of the quiver $Q$ from the
Theorem as a matrix algebra and determine its admissible ideals.  For
the moment, $k$ can be any field.  Regarding the multiplication of
arrows we follow the convention from \cite{ARS}: With respect
to the primitive idempotents $e_{i}$, corresponding to the vertices
${\bf i}$ for $i=1,2,3$, the path algebra $kQ = k\langle
e_{1},e_{2},e_{3};\alpha,\beta,\gamma\rangle$ is defined by the
relations
$$
        e_{1}\alpha =\alpha = \alpha e_{2}\;, \quad
        e_{2}\beta = \beta = \beta e_{2}\;, \quad
        e_{3}\gamma =\gamma = \gamma e_{2}\;,\quad
        0 =\gamma\alpha = \beta \alpha = \gamma\beta\ .
$$

Throughout, $k[\beta]$ denotes the polynomial algebra in one
variable $\beta$ over $k$, and $(\beta^n)$ 
denotes the ideal generated in $k[\beta]$ by $\beta^n$.
 
\begin{lemma}
        The path algebra
        $kQ$
        can be realized as a triangular matrix algebra
        $$
        \vp:kQ\xrightarrow{\cong}
        \left(\begin{matrix}
        k& k[\beta] & k[\beta]\\
        0&k[\beta]&k[\beta]\\
        0&0&k
        \end{matrix}\right)
        $$
        where $\vp$ is uniquely determined through
        \begin{align*}
                \vp(xe_{1} + y e_{2} + ze_{3}+ u\alpha + v \beta + w\gamma)=
                \left(\begin{matrix}
                x& u & 0\\
                0&y + v\beta&w\\
                0&0&z
                \end{matrix}\right)
        \end{align*}
        for $x,y,z,u,v,w\in k$.
    Under this isomorphism the powers of the path ideal
    $(\alpha,\beta,\gamma)$ are mapped to
    \begin{align*}
                \vp\left((\alpha,\beta,\gamma)^{n}\right)=
                \left(\begin{matrix}
                0& (\beta^{n-1}) & (\beta^{\max(0,n-2)})\\
                0&(\beta^{n})&(\beta^{n-1})\\
                0&0&0
                \end{matrix}\right) \quad \text{for}\quad n>0\,.
        \end{align*}
\qed
\end{lemma}
An ideal $I$ in $kQ$ is {\em admissible\/} if
$(\alpha,\beta,\gamma)^{N}\subseteq I \subseteq
(\alpha,\beta,\gamma)^{2}$ for some $N\ge 2$. If $k$ is
algebraically closed, the quotients $kQ/I$ of $kQ$  by
admissible ideals $I$ are up to Morita equivalence precisely 
the finite dimensional $k$-algebras
with that ordinary quiver. Hochschild (co-)homology is Morita 
invariant, see \cite[(1.2.4), (1.5.6)]{Lod}, a fact established for
outer derivations of matrix algebras already by P.~Dirac \cite{Di} in 1925!

\begin{lemma}
        The admissible ideals in $kQ$ are precisely
                $$
                I(n;n',n'';V)=
                (\alpha\beta^{n'},\beta^{n},\beta^{n''}\gamma, \alpha V\gamma)
                \subset kQ
                $$
                where $2\le n\ ;\  1\le n',n''\le n$ and $V\subseteq
k[\beta]$ is a
                $k$-subspace that contains $(\beta^{m})$ with
$m = \min(n',n'')$.
                Such an ideal is mapped to
                $$
                \vp\left(I(n;n',n'';V) \right) =
                \left(\begin{matrix}
                        0& (\beta^{n'}) & V\\
                        0&(\beta^{n})&(\beta^{n''})\\
                        0&0&0
                        \end{matrix}\right)\ .
                $$
                If $\artA =kQ/I(n; n'; n''; V)$ is the corresponding 
                algebra, then its {\em Cartan matrix\/} is
                \begin{equation*}
                         (\dim_{k} e_{i}\artA e_{j})_{1\le i,j\le 3} =
                         \left(\begin{matrix} 1
                         & n' & m-d\\
                         0 & n & n''\\
                         0&0&1
                \end{matrix}\right)
       \end{equation*}
       where $d=\dim _kV/(\beta^{m})$.
       \qed
\end{lemma}
The vectorspace $V$ is uniquely determined by its image
$V/(\beta^{m})\subseteq k[\beta]/(\beta^{m})$,
equivalently, by the corresponding point
$$
        \left(\frac{V}{(\beta^{m})}\subseteq
        \frac{k[\beta]}{(\beta^{m})}\right)\in \Grass_{k}(d,m)
$$
in the Grassmanian of vector subspaces of dimension
$d$ of the $m$-dimensional $k$-vectorspace
$k[\beta]/(\beta^{m})$.  This {\em continuous invariant\/}
ranges hence over the $k$-rational points of an irreducible projective
algebraic variety of dimension $d(m-d)$.

\begin{sit}
\label{Cartan}
	The determinant of the Cartan matrix being equal to $n\ge 2$, this
	matrix is not invertible over the integers, whence the global
	dimension of $A$ is necessarily infinite.  This fact is as well
	implied immediately by the more general affirmative solution to the
	{\em ``no loops conjecture''} due to H.~Lenzing \cite{Len} (see also
	K.~Igusa \cite{Ig}).
\end{sit}

\begin{remark}
        The Hochschild homology of $\artA =kQ/I(n;n',n'';V)$ over $k$ depends
        solely upon its Loewy length; as $k$-vectorspaces,
                $$
                \Hoch_{*}(\artA ) \cong
\Hoch_{*}(k)\oplus\Hoch_{*}\left(k[\beta]/(\beta^{n})\right)\oplus
                \Hoch_{*}(k)\,.
                $$
		This follows for example from applying Theorem 1.2.15 in 
		\cite{Lod} twice.  In particular, $\Hoch_{i}(\artA ) \cong 
		\Hoch_{i}\left(k[\beta]/(\beta^{n})\right)$ for $i > 0$ and 
		$\Hoch_{i}\left(k[\beta]/(\beta^{n})\right)\ne 0$ for any $i$ 
		as soon as $n \ge 2$; \cite[Prop. 5.4.15]{Lod} 
		gives an explicit description of these groups.
\end{remark}

\begin{remark}
	Note also that for any field extension $k'$ of $k$ and any finite 
	dimensional $k$-algebra $\artA$ one has a natural $k'$-linear 
	isomorphism 
	$\Hoch^{*}(\artA/k)\otimes_{k}k'\cong 
	\Hoch^{*}\left((\artA\otimes_{k}k')/k'\right)$
	whence vanishing can be tested over any field containing $k$.
\end{remark}

\section{The Derivations}
Now we investigate derivations.  First note that for a
finite dimensional basic algebra $\artA$ over $k$, modulo inner
derivations
every derivation
vanishes on the primitive idempotents.  A derivation that vanishes on
these idempotents is {\em normalized\/}; it respects the two-sided
Pierce decomposition of $\artA$.  In particular,
$\Hoch^{1}(\artA)$ is isomorphic to the vectorspace of normalized
derivations modulo the inner derivations  vanishing on the
primitive idempotents, see \cite[(3.1)]{Hap}. Recall as well that
$\Hoch^{1}(\artA)$ inherits the Lie
$k$-algebra structure from the Lie algebra of all (normalized)
derivations.

As the arrows together with the primitive idempotents generate the
path algebra $P$ of any quiver, a normalized derivation $D$ of such
algebra has a unique representation
$$
D =\sum_{\alpha:{\bf i}\to {\bf j}}x_{\alpha}
\frac{\partial}{\partial \alpha}
$$
where $x_{\alpha}:=D(\alpha)$ lies in the same two-sided Pierce 
component as the arrow $\alpha$.  If $J\subset P$ is an ideal, then 
$D$ descends to a derivation of $P/J$ if and only if $D(J)\subseteq 
J$.

To take care of the possibly finite characteristic of $k$,
we introduce two secondary invariants of $\artA = kQ/I(n;n',n'';V)$:
$$
\delta = \gcd(n,n',n'') \quad\text{and}\quad c = \dim_{k}(k/\delta k) =
\begin{cases}
        0& \text{if $\delta\ne 0$ in $k$}\,,\\
        1& \text{if $\delta = 0$ in $k$}\,.
\end{cases}
$$
\begin{lemma}
\label{derivations}
With notation as introduced above, one has:
   \begin{enumerate}
       \item
	   \label{normalized}
	   The normalized $k$-derivations of the path algebra $kQ$ are of 
	   the form 
	   $$
       D = \alpha a(\beta)\frac{\partial}{\partial \alpha} +
           b(\beta)\frac{\partial}{\partial \beta} +
           c(\beta)\gamma\frac{\partial}{\partial \gamma}
       $$
       with $a(\beta), b(\beta), c(\beta)\in k[\beta]$.

       \item
       \label{induced}
	   The normalized derivation $D$ of $kQ$ descends to a derivation on $\artA$ 
	   if and only if $\delta b(0) =0$ and the first order 
	   differential operator
       $$
	   \Theta_{D}=(a(\beta)+c(\beta)) + 
	   b(\beta)\frac{\partial}{\partial \beta}
       $$
       on $k[\beta]$ maps $V$ to itself.

       \item
       \label{inner}
	   The normalized derivation $D$ induces an inner derivation on 
	   $\artA$ if and only if
       \begin{align*}
			a(\beta)+c(\beta)\equiv a(0)+c(0)\bmod 
			\beta^{\min(n',n'')}\quad\text{and}\quad 
			b(\beta)\equiv 0\bmod \beta^{n}\,.
       \end{align*}
       \end{enumerate}
\end{lemma}

\begin{proof}
   In $kQ$, the twosided Pierce component of $\alpha$ is $e_{1}(kQ) 
   e_{2} = \alpha k[\beta]$, that of $\beta$ is $e_{2}(kQ) e_{2}= 
   k[\beta]$ and that of $\gamma$ is $e_{2}(kQ) e_{3}=k[\beta]\gamma 
   $, whence (\ref{normalized}).

   We now check what it means for $D$ to map (each twosided Pierce 
   component of) $I=I(n;n',n'';V)$ to itself: As $e_{2}I e_{2}= 
   (\beta^{n})$, we get
   $$
   D(\beta^{n}) = n b(\beta)\beta^{n-1}\in (\beta^{n}) \quad \text{if 
   and only if} \quad n b(0) = 0\,.
   $$
   Applying $D$ to $\alpha\beta^{n'}$ yields
   $$
   D(\alpha\beta^{n'})= \alpha \left(a(\beta)\beta^{n'} +
   n' b(\beta)\beta^{n'-1}\right)\ ,
   $$
   which is in $\alpha (\beta^{n'})= e_{1}Ie_{2}$ if and only if 
   $n'b(0)=0$.  Similarly applying $D$ to $\beta^{n''}\gamma$ yields 
   the condition $n''b(0) = 0$.  Now $nb(0)=n'b(0)=n''b(0)=0$ if and 
   only if $\delta b(0)=0$.  Finally, for $\alpha v\gamma$ with $v\in 
   V\subseteq k[\beta]$, one has
   $$
   D(\alpha v \gamma)= \alpha\left(a(\beta) v + b(\beta)\frac{\partial 
   v}{\partial \beta} + c(\beta) v\right )\gamma = \alpha 
   \Theta_{D}(v) \gamma
   $$
   and (\ref{induced}) follows.

   If $a = a_{11}e_{1} + a_{22}e_{2} + a_{33}e_{3} + \alpha a_{12} + 
   a_{23}\gamma + \alpha a_{13}\gamma \in kQ$ with $a_{11}, a_{33}\in 
   k$ and the remaining coefficients in $k[\beta]$, then 
   the inner derivation $[a,\ ]$ vanishes on the primitive idempotents 
   if and only if $\alpha a_{12}=a_{23}\gamma =\alpha a_{13}\gamma =0$ 
   in $\artA $.  
   Calculating its value on the arrows yields 
   $$[a,\ ] = 
   \alpha(a_{11}- a_{22})\frac{\partial}{\partial \alpha} + (a_{22}- 
   a_{33})\gamma\frac{\partial}{\partial \gamma}
   $$ 
   and (\ref{inner}) follows.
\end{proof}

To summarize this detailed information succinctly, we interpret it in
terms of differential operators of first order on $k[\beta]$.  
If $p$ is a polynomial in $\beta$, its derivative with
respect to $\beta$ is denoted by $p'$ as usual.

\begin{lemma}
   Consider
   $
   J_{1}\subseteq J_{2}\subseteq V \subseteq k[\beta] 
   $
   where $J_{\nu}=(\beta^{e_{\nu}}), \nu=1, 2$, are proper ideals and 
   $V$ is some vector subspace in $k[\beta]$.

   If $e\in \bbbn$ divides $e_{\nu}$ for $\nu=1,2$, then the 
   $k$-vectorspace
   $$
   \cald_{V,e}=\left\{\left. \Theta= A(\beta) + 
   B(\beta)\dis\frac{\partial}{\partial \beta} \right|\begin{matrix} 
   A,B\in k[\beta], e B(0) =0, \Theta(V)\subseteq V
   \end{matrix}
   \right\}
   $$
   is a Lie subalgebra of the Lie $k$-algebra of all first order 
   differential operators on $k[\beta]$ that contains $J=\left( 
   k\oplus J_{2}\right)\bigoplus J_{1}\dis\frac{\partial}{\partial 
   \beta}$ as an ideal and $\cald_{V,1}$ as a Lie subalgebra.

   The quotient Lie $k$-algebra $\cald(V,e;J_{1},J_{2})=\cald_{V,e}/J$ 
   is finite dimensional and contains an ideal isomorphic to 
   $\left(J_{2}/J_{1}\right)\dis\frac{\partial}{\partial \beta}$.  Its 
   dimension satisfies
   $$
   0\le e_{1}-e_{2}\le \dim_{k} \cald(V,e;J_{1},J_{2}) - \epsilon \le 
   e_{1}+e_{2}- 2
   $$
   where $\epsilon=\dim_{k}(k/e k)$.
\end{lemma}

\begin{proof}
   It suffices to recall that the Lie bracket of two first order 
   differential operators on $k[\beta]$ is given by
   $$
   \left[A_{1}(\beta) + B_{1}(\beta)\dis\frac{\partial}{\partial 
   \beta}, A_{2}(\beta) + B_{2}(\beta)\dis\frac{\partial}{\partial 
   \beta}\right]= \det\left(
        \begin{matrix}
                B_{1}& A'_{1}\\
                B_{2}& A'_{2}\\
        \end{matrix}\right)
        +
        \det\left(
        \begin{matrix}
                B_{1}& B'_{1}\\
                B_{2}& B'_{2}\\
        \end{matrix}
        \right)\dis\frac{\partial}{\partial \beta}\,.
   $$
   Now the ideal $J_{\nu}=(\beta^{e_{\nu}})$ is mapped to itself by 
   the vectorfield $B(\beta)\dis\frac{\partial}{\partial \beta}$ if 
   and only if $e_{\nu}B(0)=0$.  The final assertion follows from the 
   inclusions
   $$
   \left( k\oplus J_{2}\right)\bigoplus
   J_{1}\dis\frac{\partial}{\partial \beta}\subseteq
   \left( k\oplus J_{2}\right)\bigoplus
   J_{2}\dis\frac{\partial}{\partial \beta}\subseteq
   \cald_{V,1}\subseteq
   \cald_{V,e}\,.
   $$
\end{proof}

For a finite dimensional algebra $\artA = kQ/I(n;n',n'';V)$ as before
we have thus the following description of $\Hoch^{1}(\artA )$ where we 
set again $m=\min(n', n''),\, \delta = \gcd(n,n',n'')$ and 
$c=\dim_{k}(k/\delta k)$.
\begin{corollary}
	Associating to a normalized derivation $D$ on $kQ$ the first order 
	differential operator $\Theta_{D}$ on $k[\beta]$ as in {\em 
	(\ref{derivations}.\ref{induced})\/} induces an isomorphism of Lie 
	$k$-algebras
    $$
    \Hoch^{1}(\artA )\cong \cald(V,\delta; \beta^{n},\beta^{m}) \,.
    $$
	In particular, $\cald'(\artA) = \cald(V,1; \beta^{m},\beta^{m})$ 
	is a Lie algebra subquotient of $\Hoch^{1}(\artA)$ and
    $$
    \dim_{k}\Hoch^{1}(\artA) = \dim \cald'(\artA ) + n - m + c\le (m-1) +
    (n-1) + c\,.
    $$
\qed
\end{corollary}

\begin{example}
\label{ideal}
	If $ V$ is an {\em ideal\/} in $k[\beta]$, so that $V = 
	(\beta^{m-d})$, and if $\delta$ divides $d$, then the upper bound 
	is achieved,
    $$
        \dim \Hoch^{1}(\artA) = (n-1) + (m - 1) + c> 0\,.
    $$
    Note that the hypothesis is trivially satisfied if $m=1$, in which
    case $c=0$ and $\dim \Hoch^{1}(\artA) = n-1$.
\end{example}

\begin{sit}
\label{D'}
	To investigate the structure of $\cald'(\artA )= \cald(V,1; 
	\beta^{m},\beta^{m})$, choose polynomials $p_{1}, p_{2},\ldots , 
	p_{d}\in k[\beta]$ that generate the vectorspace $V$ minimally modulo 
	$(\beta^{m})$.  As the classes of the polynomials $p_{i}$ in 
	$k[\beta]/(\beta^{m})$ are $k$-linearly independent, $\cald'(\artA 
	)$ is isomorphic as vectorspace to the solutions of the system of 
	$d$ equations
\begin{equation}
\label{thesystem}
	\left(\begin{matrix} 
	p_{1}& p'_{1}\\
    p_{2}& p'_{2}\\            
    \vdots&\vdots\\
    p_{d}& p'_{d}\\
    \end{matrix}\right)
    \left(\begin{matrix}
    A\\
	B 
	\end{matrix}\right) 
	\equiv 
	\left(\begin{matrix} 
	a_{11}&\ldots & a_{1d}\\
    a_{21}&\ldots & a_{2d}\\
    \vdots&\ddots&\vdots\\
    a_{d1}&\ldots & a_{dd}
    \end{matrix}\right)
    \left(\begin{matrix}
    p_{1}\\
    p_{2}\\
    \vdots\\
    p_{d}
    \end{matrix}\right)\bmod \beta^{m}
\end{equation}
with $a_{ij}\in k$ and $A,B\in k[\beta]; A(0)=B(0)=0; \deg A, \deg B
\le m-1$.

Comparing in each of these equations the coefficients of $\beta^{i}$
for $i=0,\ldots,m-1$ yields $m$ linear equations over $k$; thus we
obtain in total $dm$ equations in $d^{2}+ 2(m-1)$ unknowns, whence
the dimension of $\cald'(\artA )$ satisfies
\begin{equation}
\label{corank}
        \dim \cald'(\artA ) \ge d^{2}+ 2(m-1) - dm = (d-2)(d+2-m) + 2 \,.
\end{equation}
Evaluating this lower bound yields immediately the following.
\begin{corollary}
\label{linearbound}
	Whenever $d \le 2\le m$ or $d > 2$ and $d + \dis\frac{2}{d-2} > m-2$, 
	then $\cald'(\artA )\ne 0$ and so $\dim_{k}\Hoch^{1}(\artA )> n-m 
	+ c \ge 0$.
\end{corollary}

\begin{proof}
	If $d=0$, then the right hand side of (\ref{corank}) evaluates to 
	$2 m - 2 \ge 2$, if $d=1$ it evaluates to $m-1\ge 1$, if $d=2$ it 
	evaluates to $2$. The second case simply rewrites $(d-2)(d+2-m) + 
	2>0$.
 \end{proof}

Thus to find examples with $\Hoch^{1}(\artA )=0$, one needs $d\ge 3$
and $m=n\ge 7$. Note that in this case the algebra $\artA=kQ/I(n;n,n;V)$
is of infinite representation type; see for example
the classification of the maximal algebras with 2 simple modules in 
\cite{BG}.
So we get the following.

\begin{corollary}
\label{frt}
Let $\artA$ be a finite dimensional $k$-algebra with ordinary quiver
$Q$ as in Theorem {\em (\ref{main result})\/}. If
$\artA$ is of finite representation type, then $\Hoch^{1}(\artA)$ does not
vanish.
\end{corollary}

\end{sit}

\section{Examples in Large Characteristic}

If the field $k$ has a large enough prime field, we can use ``finite 
Fourier analysis'' to find examples where $\cald'(\artA)=0$.  For 
$y\in k$, let $e^{y\beta}\bmod \beta^{m}$ denote the image of 
$$
e^{y\beta}=\sum_{\nu}y^{\nu}\dis\frac{\beta^{\nu}}{\nu!}\in
k[[\beta]]
$$
in $k[[\beta]]/(\beta^{m})\cong k[\beta]/(\beta^{m})$ for $m\in 
\bbbn$.  This class is well defined as soon as $(m-1)!\ne 0$ in $k$.  
Clearly, if defined, $e^{y\beta}\bmod \beta^{m}$ is a unit in 
$k[\beta]/(\beta^{m})$ with inverse $e^{-y\beta}\bmod \beta^{m}$.

\begin{lemma}
\label{Vandermonde}
        If $(m-1)!\ne 0$ in $k$ and $Y\subseteq k$ is a
        subset of cardinality $\mu$, then the classes $\{e^{y\beta}\bmod
        \beta^{m}\}_{y\in Y}$
        span a vectorspace of dimension $d=\min(m,\mu)$ in
$k[\beta]/(\beta^{m})$.
\end{lemma}

\begin{proof}
        Vandermonde's determinant.
\end{proof}

\begin{proposition}
\label{exponential example}
        Let $k$ be a field and $d, m$ integers with $d\ge 3, m \ge 3d-2$,
        and $(m-1)!\ne 0$ in $k$.  Assume $Y\subseteq k$ is a finite
        subset of cardinality $d$ such that the set of differences
        $$
        Y-Y = \{y_{1}-y_{2}\mid y_{1},y_{2}\in Y\}\subseteq k
        $$
        has maximal cardinality, equal to $1 + 2 \binom{d}{2}$.
        Let $V$ be the $k$-vectorspace given by
        $(\beta^{m})\subseteq V\subseteq k[\beta]$ and
        $$
        V/(\beta^{m})= \sum_{y\in Y} k\left( e^{y\beta}\bmod
        \beta^{m}\right)\subseteq k[\beta]/(\beta^{m})\,.
        $$
        There is then no nonzero differential operator $\Theta=A(\beta) +
        B(\beta)\dis\frac{\partial}{\partial \beta}$ with $A,B\in
        (\beta)$ and\ $\deg A,\deg B< m$ that transforms the
        $k$-vectorspace $V$ into itself.
\end{proposition}

\begin{proof}
        The condition on (the characteristic of) $k$ ensures that $V$ is
        well defined.  If $\Theta=A(\beta) +
        B(\beta)\dis\frac{\partial}{\partial \beta}$ is a
        differential operator on $k[\beta]$ with $B(0)=0$, then
    \begin{equation}
    \label{exponential}
            \Theta(e^{y\beta}\bmod \beta^{m}) = \left(A(\beta) + y
            B(\beta)\right)e^{y\beta}\bmod \beta^{m}\,.
    \end{equation}
        Now assume
        \begin{equation}
        \label{expsystem}
                \left(A(\beta) + y
                B(\beta)\right)e^{y\beta}\bmod \beta^{m} =\sum_{y'\in
                Y}a_{yy'} e^{y'\beta}\bmod \beta^{m}
        \end{equation}
        for some matrix $(a_{yy'})_{y,y'\in Y}$ of elements from $k$.
        Multiplying by $e^{-y\beta}\bmod
        \beta^{m}$ and subtracting yields an equation
        \begin{equation}
        \label{equation for B}
                \left(y_{1}-y_{2}\right)B = \sum_{y'\in Y}a_{y_{1}y'}
                e^{(y'-y_{1})\beta} -\sum_{y'\in Y}a_{y_{2}y'}
                e^{(y'-y_{2})\beta}\bmod \beta^{m}
        \end{equation}
        for each pair of elements $(y_{1},y_{2})\in Y\times Y$.  If $y_{3}$ is
        a third element from $Y$, there result equations
        \begin{align*}
                &\left(y_{1}-y_{2}\right)\left(\sum_{y'\in Y}a_{y_{1}y'}
                e^{(y'-y_{1})\beta} -\sum_{y'\in Y}a_{y_{3}y'} 
                e^{(y'-y_{3})\beta}\right)\bmod \beta^{m} \\
                =&\left(y_{1}-y_{3}\right)\left(\sum_{y'\in Y}a_{y_{1}y'}
                e^{(y'-y_{1})\beta} -\sum_{y'\in Y}a_{y_{2}y'}
                e^{(y'-y_{2})\beta}\right)\bmod \beta^{m}
        \end{align*}
        in $k[\beta]/(\beta^{m})$.  The (classes of) exponential
        functions involved are $1 = e^{0\beta}$ and $e^{y_{i}-y_{j}}$ for
        $y_{i}\in Y\setminus\{y_{j}\}$ with $j=1,2,3$.  By the assumption on
        the set of differences, if $y_{1},y_{2},y_{3}$ are pairwise
        distinct, the
        set
        $$
        Z:= \{0, \, y_{i}-y_{j};\, y_{i}\in
        Y\setminus\{y_{j}\}; \;\; j=1,2,3\}\subseteq Y - Y\subseteq k
        $$
        has cardinality equal to $1 + 3(d-1) = 3d-2$, and as $m\ge 3d-2$,
        the corresponding classes $\{e^{z\beta}\bmod
        \beta^{m}\}_{z\in
        Z}$ are linearly independent in $k[\beta]/(\beta^{m})$ by Lemma
        (\ref{Vandermonde}).  It
        follows in particular that $a_{y_{i}y_{2}}=0$ for $y_{i}\ne
        y_{2}$.  As the choice of $y_{2}$ was arbitrary and as there are at
        least $3$ distinct elements in $Y$, it follows that $a_{yy'}=0$
        whenever $y\ne y'$. In turn,
        the system (\ref{expsystem}) evaluates now at $\beta=0$ to
        $$
        A(0) + y B(0) = \sum_{y'\in Y}a_{yy'} = a_{yy}\;,\quad y\in Y,
        $$
        and so $A(0) = B(0) = 0$ implies that the diagonal elements $a_{yy}$
        vanish as well. Finally, equation (\ref{equation for B}) shows
        that $B=0$ and any equation in (\ref{expsystem}) yields $A=0$.
\end{proof}

\begin{corollary}
\label{exponential examples}
        Let $I\subset kQ$ be an admissible ideal with discrete invariants
        $(n,n',n'',d)$ and continuous invariant
        $\left(V/(\beta^{m})\right)\in \Grass_{k}(d,m)$.  If $m!\ne
        0$ in $k$ for $m=\min(n',n'')$, and if $m\ge 3d-2 \ge 7$, then
        there exists a Zariski open and dense subset ${\mathfrak U}\subset
        \Grass_{k}(d,m)$ such that $\artA = kQ/I$ satisfies
        $$
        \dim_{k}\Hoch^{1}(\artA) = n-m
        $$
        whenever $\left(V/(\beta^{m})\right)\in\mathfrak U$.
\end{corollary}

\begin{proof}
	The conditions guarantee that $d,m$ satisfy the hypotheses of 
	Proposition (\ref{exponential example}) and that $\delta = 
	\gcd(n,n',n'')\le m$ is a unit in $k$, whence $c=0$.  As 
	$\dim_{k}\Hoch^{1}(\artA)$ is upper semicontinuous on 
	$\Grass_{k}(d,m)$ it remains to show that the minimal value is 
	taken on generically, for example over the purely transcendental 
	field extension $k' = k(y_{1},\ldots,y_{d})$ of $k$.  There 
	Proposition (\ref{exponential example}) applies.
\end{proof}

\begin{example}
\label{explicit exponential}
	The simplest case occurs for $d=3, m=7$.  With $Y=\{0,1,3\}$, the 
	set of differences $Y-Y = \{-3,-2,-1,0,1,2,3\}$ has maximal 
	cardinality and so for any field $k$ of characteristic greater 
	than $7$ or of characteristic zero, the $k$-algebra
    $$
	\artA = kQ/(\beta^{7}, \alpha\gamma , \alpha (e^{\beta}\bmod 
	\beta^{7})\gamma , \alpha (e^{3\beta}\bmod \beta^{7})\gamma)
    $$
	of dimension $27$ has a loop in its ordinary quiver but satisfies 
	$\Hoch^{1}(\artA) = 0$.
\end{example}

\section{Arbitrary Characteristic}

To find examples in arbitrary characteristic, we analyze a bit further 
the system of equations (\ref{thesystem}). One may choose a $k$-basis 
$\left(p_{i}(\beta)\right)_{i\in \bbbn}$ of $V\subseteq k[\beta]$ such 
that the {\em orders\/} $\ord p_{i}$ of these polynomials satisfy
$$
0\le \nu_{1}=\ord p_{1} < \nu_{2}=\ord p_{2} < \dotsb < \nu_{i}=\ord
p_{i} <\dotsb\,.
$$
As $(\beta^{m})\subseteq V$ one has $\nu_{i}\le m+i-1$; the classes of 
$p_{1},\ldots,p_{d}$ constitute a basis of $V/\left(\beta^{m}\right)$.  
In view of Example (\ref{ideal}) we assume henceforth that $m \ge 2$.

\begin{lemma}
	With the preceding notations and assumptions, the dimension of the 
	Lie algebra $\cald'(\artA )$ from {\em (\ref{D'})} satisfies
    $$
	\dim \cald'(\artA ) \ge 2m - \max(1, \nu_{1} ,m-\nu_{2} + 1) - 
	\max(1, m-\nu_{1})\,.
    $$
    In particular, $\cald'(\artA ) = 0$ implies $\nu_{1}=0$ and $\nu_{2}=1$.
\end{lemma}

\begin{proof}
	Consider first the integers $\mu$ with $m\ge \mu > \max(1,\nu_{1}, 
	m-\nu_{2} + 1)$.  As $\mu > \nu_{1}$, each $B_{\mu}(\beta) = 
	-\beta^{\mu-\nu_{1}-1} p_{1}$ is a nonzero polynomial.  As $m\ge 
	\mu > 1$ the order satisfies
    $$
    m > \ord B_{\mu} = (\mu-\nu_{1}-1) + \nu_{1} = \mu -1 > 0\,
    $$
	whence $B_{\mu}(0)=0$ and the classes $\{B_{\mu}(\beta) \bmod 
	{\beta^{m}}\}_{\mu}$ are $k$-linearly independent in 
	$k[\beta]/(\beta^{m})$.  For the differential operator 
	$\Theta_{\mu} = \beta^{\mu-\nu_{1}-1}\left(p'_{1} - 
	p_{1}\frac{\partial}{\partial\beta}\right)$ one has 
	$\Theta_{\mu}(p_{1})=0$ by construction and for $i\ge 1$ one finds
    \begin{align*}
		\ord(\Theta_{\mu}(p_{i})) &= 
		\ord\left(\beta^{\mu-\nu_{1}-1}\left(p'_{1} p_{i} - 
		p_{1}p'_{i}\right)\right)\\
		& \ge (\mu-\nu_{1}-1) + (\nu_{1}+ \nu_{i} - 1)\ge \mu + 
		\nu_{2} - 2 \ge m
    \end{align*}
	in view of $\nu_{i}\ge \nu_{2}$ and $\mu \ge m-\nu_{2} + 2$.  It 
	follows that $\Theta_{\mu}(p_{i})\in (\beta^{m})\subseteq V$ for 
	each $i$, whence the subset $\{\Theta_{\mu}\}_{\mu}\subseteq 
	\cald'(\artA )$ is $k$-linearly independent of cardinality $m - 
	\max(1, \nu_{1} ,m-\nu_{2} + 1)$.

	Now consider the integers $\mu'$ with $m >\mu' \ge \max(1, 
	m-\nu_{1})$.  The polynomials $A(\beta) = \beta^{\mu'}$ satisfy 
	then $A(0) = 0$ and $\ord \left(Ap_{i}\right) = \mu' + \nu_{i}\ge 
	\mu' + \nu_{1} \ge m$, whence the set of (scalar) differential 
	operators $\{\Theta^{\mu'} = \beta^{\mu'}\}_{\mu'}$ defines a 
	$k$-linearly independent subset of $\cald'(\artA )$ of cardinality 
	$m - \max(1, m-\nu_{1})$.  As the union 
	$\{\Theta_{\mu}\}_{\mu}\cup \{\Theta^{\mu'}\}_{\mu'}$ is clearly 
	still linearly independent over $k$, the lower bound on $\dim 
	\cald'(\artA )$ follows.

	For the final assertion observe that $0\le \nu_{1}\le m$ whence 
	$1\le \max(1,m-\nu_{1})\le m$, with equality on the right if and 
	only if $\nu_{1} =0$.  Similarly $1\le \nu_{2}\le m+1$ and so 
	$1\le \max(1,\nu_{1},m-\nu_{2}+1)\le m$ with equality on the right 
	if and only if $\nu_{1}=m$ or $\nu_{2}= 1$.  The claim follows.
\end{proof}

\begin{sit}
\label{triangular}
	Due to the preceding result we assume henceforth that 
	$\nu_{1}=0, \nu_{2}=1$. We normalize further the polynomials $p_{i}$ so 
	that
    $$
    p_{i}(\beta) = \beta^{\nu_{i}} +
        \sum_{j=1}^{m-1-\nu_{i}}p_{i,j}\beta^{\nu_{i}+j}\quad,\quad
        p_{i,j}\in k\,.
    $$
	As $\nu_{2}=1$, we may moreover assume $p_{1,1}=0$.
	Indeed, if the {\em order sequence\/} $(0 <1 <\nu_{3} 
	<\ldots <\nu_{d})$ is fixed, we can impose $p_{i,j}= 
	0$ whenever $j = \nu_{i'}-\nu_{i}$ for some $i'>i$, and in that 
	case the undetermined $p_{i,j}$ serve as $d(m-d)$ affine coordinates on 
	$\Grass_{k}(d,m)$ for the classes of vectorspaces with the given 
	order sequence.
\end{sit}	
	
\begin{sit}
	Just comparing orders 
	and lowest order terms of the polynomials on both sides of the 
	equations in (\ref{thesystem}) shows that $A(0)=B(0)=0$ 
	implies
    \begin{enumerate}
    \item[(a)] $a_{ij}=0$ for $1\le j < i \le d$,

    \item[(b)] $a_{11}= 0$,

    \item[(c)] $B = a_{22}\beta + O(\beta^{2})$, and

    \item[(d)] $a_{ii} = \nu_{i} a_{22}$ for $i = 3,\ldots,d$.
    \end{enumerate}
\end{sit}

\begin{sit}
	Now consider the polynomials
    \begin{equation}
    \label{Wronski}            
       \begin{split}
       \Delta_{ij}&=\det\left(\begin{matrix}
                    p_{i}& p'_{i}\\
                    p_{j}& p'_{j} \end{matrix}\right)\\
				  &= \nu_{ji}\beta^{\nu_{j}+\nu_{i}-1}+ 
				  \left[\left(\nu_{ji}+1\right)p_{j,1}+ 
				  \left(\nu_{ji}-1\right)p_{i,1}\right]\beta^{\nu_{j}+\nu_{i}} 
				  +O(\beta^{\nu_{j}+\nu_{i}+1})
        \end{split}
   \end{equation}
   in $k[\beta]$ where $\nu_{ji}=\nu_{j}-\nu_{i}$.  The assumption 
   $\ord p_{1} = 0, \ord p_{2} = 1$ yields that $\Delta_{12}$ is a 
   unit in $k[\beta]/(\beta^m)$.  Moreover
   \begin{equation}
   \label{Pluecker}
       \left(\Delta_{2i},-\Delta_{1i},\Delta_{12}\right)
       \left(\begin{matrix}
                p_{1}& p'_{1}\\
                p_{2}& p'_{2}\\
                p_{i}& p'_{i}\\
	   \end{matrix}\right) =
	   \left(0,0\right)\,.
   \end{equation}
   Multiplying the given system (\ref{thesystem}) from the left by the 
   $d\times d$-matrix of polynomials
   \begin{equation*}
	   \left(
	   \begin{matrix}
           p'_{2}&-p'_{1}& 0  &\ldots &0  &\ldots &0  \\
           -p_{2}&  p_{1}& 0  &\ldots &0 &\ldots &0 \\
	   \Delta_{23}&-\Delta_{13}&\Delta_{12}&\ldots&0 &\ldots&0\\
	   \vdots &\vdots & \vdots &\ddots &\vdots &&\vdots\\
       \Delta_{2i}&-\Delta_{1i}&0&\ldots&\Delta_{12}&\ldots&0\\
	   \vdots & \vdots &\vdots & &\vdots &\ddots&\vdots\\
       \Delta_{2d}&-\Delta_{1d}&0&\ldots&0&\ldots&\Delta_{12}
       \end{matrix}
       \right)
   \end{equation*}
   produces the equivalent system
   \begin{equation}
   \label{Deltasystem}
       \begin{split}
       \Delta_{12}A &\equiv \sum_{j=1}^{d}\left(a_{1j}p'_{2} -
       a_{2j}p'_{1}\right)p_{j} \bmod \beta^{m}\\
       \Delta_{12}B &\equiv \sum_{j=1}^{d}\left(-a_{1j}p_{2} +
       a_{2j}p_{1}\right)p_{j} \bmod \beta^{m}\\
       0&\equiv \sum_{j=1}^{d}\left(a_{1j}\Delta_{2i} -
       a_{2j}\Delta_{1i}+
       a_{ij}\Delta_{12}\right)p_{j} \bmod \beta^{m}\;;
        \;i=3,\ldots,d\,.
       \end{split}
       \end{equation}
Using relation (\ref{Pluecker}) to write $\Delta_{1i}p_{2}= 
\Delta_{2i}p_{1} + \Delta_{12}p_{i}$ and employing the conditions found in 
(\ref{triangular}), the last $d-2$ equations in (\ref{Deltasystem}) 
become
\begin{multline}
\label{c}
        0\equiv -a_{22}\left(\Delta_{2i}p_{1} +
        (1-\nu_{i})\Delta_{12}p_{i}\right) +\\
        \quad\sum_{j=2}^{d}a_{1j}\Delta_{2i}p_{j} -
        \sum_{j=3}^{d}a_{2j}\Delta_{1i}p_{j} +
        \sum_{j=i+1}^{d}a_{i j}\Delta_{12}p_{j}\bmod \beta^{m}\;;\;
                i=3,\ldots,d\,.
\end{multline}
Consider in particular the equation for $i=d$,
\begin{equation*}
\label{d}
0\equiv -a_{22}\left(\Delta_{2d}p_{1} +
	        (1-\nu_{d})\Delta_{12}p_{d}\right) +
	        \sum_{j=2}^{d}a_{1j}\Delta_{2d}p_{j} -
	        \sum_{j=3}^{d}a_{2j}\Delta_{1d}p_{j} \bmod \beta^{m}.
\end{equation*}
If the occurring $2(d-1)$ polynomials
\begin{equation}
	\label{qs}
	\begin{split}
	 q_{d} &:= \Delta_{2d}p_{1} + (1-\nu_{d})\Delta_{12}p_{d}\\
	q_{1dj}&:= \Delta_{1d}p_{j}\quad\quad\text{for}\quad
			                  j=3,\ldots,d,\\
	q_{2dj}&:= \Delta_{2d}p_{j}\quad\quad\text{for}\quad 
						   j=2,\ldots,d,
	\end{split}
\end{equation}
are linearly independent modulo $(\beta^{m})$, then the coefficients 
$a_{1j}, a_{2j}$ have to vanish for all $j$.  Remembering that 
$\Delta_{12}$ is a unit, looking back at system (\ref{Deltasystem}) 
shows then that $A=B=0$, and so $\cald'(A)=0$.

If $d=3$, this condition is clearly as well necessary,
whence we get the following result.
\end{sit}

\begin{corollary}
\label{d=3}
	If $d=3$ and $V$ is generated modulo $(\beta^{n})$ by $3$ 
	polynomials of the form $p_{1} = 1 + O(\beta)$, $p_{2} = \beta + 
	O(\beta^{2})$ and $p_{3} = \beta^{\nu_{3}} +O(\beta^{\nu_{3}+1})$, 
	then $\artA = kQ/I(n;n,n;V)$ satisfies $\Hoch^{1}(\artA)=0$ if and 
	only if $n\ne 0$ in $k$ and the classes of the four polynomials
    \begin{align*}
         q_{3} &= \Delta_{23}p_{1} + 
         (1-\nu_{3})\Delta_{12}p_{3},\quad
        q_{133} = \Delta_{13}p_{3},\quad
        q_{232} = \Delta_{23}p_{2},\quad
        q_{233} = \Delta_{23}p_{3}
    \end{align*}
    are $k$-linearly independent in $k[\beta]/(\beta^{n})$.
    \qed
\end{corollary}

To investigate the polynomials in (\ref{qs}), look at their low order
terms. Using (\ref{Wronski}) and assuming, as we may according to
 (\ref{triangular}), that $p_{1,1}=0$
we obtain
\begin{align*}
	q_{d} &= \left(p_{d,1} - \nu_{d}p_{2,1}\right)\beta^{\nu_{d}+1} + 
	O(\beta^{\nu_{d}+2})\\
	q_{1dj}& = \nu_{d}\beta^{\nu_{d}+\nu_{j}-1} +
	O(\beta^{\nu_{d}+\nu_{j}})\qquad \text{for}\ j=3,\ldots,d\,,\\
	q_{2dj}& = (\nu_{d}-1)\beta^{\nu_{d}+\nu_{j}}+ 
	[\nu_{d}p_{d,1}+(\nu_{d}-1)p_{j,1}+(\nu_{d}-2)p_{2,1}] 
	\beta^{\nu_{d}+\nu_{j}+1} \\
	&\quad\qquad + O(\beta^{\nu_{d}+\nu_{j}+2})\qquad \text{for}\ j=2,\ldots,d\,.
\intertext{Now assume that $\ch k = p > 0$ and that 
$\nu_{d}=\ord p_{d}\equiv 1\bmod p$. In this case the preceding 
expressions simplify to}
	q_{d} &= \left(p_{d,1} - p_{2,1}\right)\beta^{\nu_{d}+1} + 
	O(\beta^{\nu_{d}+2})\\
	q_{1dj}& = \beta^{\nu_{d}+\nu_{j}-1} +
	O(\beta^{\nu_{d}+\nu_{j}+1})\qquad \text{for}\ j=3,\ldots,d\,,\\
	q_{2dj}& =  
	\left(p_{d,1} - p_{2,1}\right)
	\beta^{\nu_{d}+\nu_{j}+1} + O(\beta^{\nu_{d}+\nu_{j}+2})\qquad \text{for}\ j=2,\ldots,d\,.
\end{align*}
This shows the following result.

\begin{corollary}
\label{positive characteristic}
	Assume $\ch k = p > 0$, $d\ge 3$ and $\nu_{d}=\ord p_{d}\equiv 
	1\bmod p$.  If $p_{d,1}- p_{2,1}\ne 0$, and if furthermore the set
	$$
	\{\nu_{d}+1,\nu_{d}+\nu_{j}-1;j=3,\ldots,d; 
	\nu_{d}+\nu_{j}+1;j=2,\ldots,d\}
	$$ 
	contains $2(d-1)$ distinct integers smaller than $n$, then $\artA= 
	kQ/(n;n,n;V)$ satisfies $\Hoch^{1}(\artA)=0$ whenever $n\ne 0$ in 
	$k$.
\end{corollary}

\begin{example}
\label{positive example}
	For $d=3$ the set in question is just $\nu_{3} + \{1,\nu_{3}-1, 2, 
	\nu_{3}+1\}$ and so the Corollary (\ref{positive characteristic}) 
	applies to the algebra $\artA =kQ/I(6p+1,6p+1,6p+1,V)$ where still 
	$\ch k = p > 0$, and $V/(\beta^{6p+1})$ is generated by $p_{1}= 
	1$, $p_{2}= \beta+\beta^{2}$, $p_{3}= \beta^{2 p+1}$.  Explicitly,
    $$
	\artA =kQ/\left(\beta^{6p+1}, \alpha\gamma, 
	\alpha(\beta+\beta^{2})\gamma, \alpha\beta^{2p+1}\gamma\right)\,.
    $$
\end{example}

To summarize, we now have the\\
{\sc Proof of Theorem} (\ref{main result}): Combine (\ref{Cartan}), 
(\ref{frt}), (\ref{explicit exponential}) and (\ref{positive 
example}).\qed

\section{Concluding Remarks}

\begin{sit}
	An unsatisfactory feature of the example in (\ref {positive example}) 
	is perhaps that its Loewy length grows with the characteristic of 
	the coefficient field.  One might thus want to analyze directly the 
	situation for $d=3$ in any characteristic. For
	\begin{align*}
		p_{1} &= 1 + p_{1,3}\beta^{3} + p_{1,4}\beta^{4} +
		p_{1,5}\beta^{5} + p_{1,6}\beta^{6} + p_{1,7}\beta^{7} +\ldots \\
		p_{2} &= \beta  + p_{2,2}\beta^{3} + p_{2,3}\beta^{4} +
		p_{2,4}\beta^{5} + p_{2,5}\beta^{6} + p_{2,6}\beta^{7}+\ldots \\
		p_{3} &= \beta^{2}  + p_{3,1}\beta^{3} + p_{3,2}\beta^{4} +
		p_{3,3}\beta^{5} + p_{3,4}\beta^{6} + p_{3,5}\beta^{7} +\ldots
	\end{align*}
	the four polynomials $q_{3}, q_{133}, q_{232}, q_{233}$ from 
	(\ref{d=3}) are of order at least $\nu_{3}+1 = 3$ and
	a quick calculation in, say,  {\sc Maple} shows that the matrix of 
	their coefficients with respect to
	$\beta^{3}, \beta^{4}, \beta^{5}, \beta^{6}$ has determinant 
	$$
	-3\left 
	(p_{1,3}-2p_{2,3}+p_{3,3}+p_{3,1}(3p_{2,2}+2p_{3,1}^2-3p_{3,2})\right)^2
	$$
	which yields the following computer aided result: There are examples 
	with vanishing first Hochschild cohomology of type $kQ/I(7;7,7;V)$ and 
	$d=3$ if and only if the characteristic of $k$ is different from $3$ 
	or $7$.  If one looks at the same coefficient matrix but with respect 
	to $\beta^{3}, \beta^{4}, \beta^{5}, \beta^{6}, \beta^{7}$, then a 
	similar computer calculation shows that the cofactor of $\beta^{5}$ is 
	a primitive polynomial in the $p_{i,j}$, whence there are examples 
	with vanishing first Hochschild cohomology of type $kQ/I(8;8,8;V)$ 
	with $d=3$ in any characteristic different from $2$.
\end{sit}

\begin{sit}
	We finish with an exercise for the diligent reader: Let $Q_{(r)}$ be 
	the quiver that has again vertices $\bf 1,2,3$ and arrows 
	$\alpha\colon{\bf 2}\to {\bf 1},\ \gamma\colon{\bf 3}\to {\bf 2}$, but 
	this time $r$ loops $\beta_{\rho}\colon{\bf 2}\to {\bf 2}$ for 
	$\rho=1,\ldots,r$.
	
	Extending the method of (\ref{exponential example}) to several 
	variables shows that in large enough or zero characteristic there are 
	finite dimensional $k$-algebras $\artA$ on $Q_{(r)}$ that satisfy 
	$\Hoch^{1}(A)=0$. Thus there is not even an upper bound on 
	the number of loops in the ordinary quiver of an Artin algebra without outer derivations.
\end{sit}


\begin{thebibliography}{ARS\ }

\bibitem{ARS}
{M.~Auslander,~I.~Reiten, and S.~Smal\o:}
{\em Representation theory of Artin algebras\/},
Cambridge Studies in Advanced Mathematics, {\bf {36}}
(University Press, Cambridge 1995).

\bibitem{BM}
{M.~J~ Bardzell and E.~N.~Marcos}, {\em Induced boundary maps for the
cohomology of monomial and Auslander algebras\/}, Canadian Mathematical
Society Conference Proceedings, {\bf 24} (1998), 47--54.

\bibitem{BG}
{K. Bongartz and P.~Gabriel}, {\em Covering Spaces in representation
Theory\/}, Invent. Math., {\bf 65} (1982), 331--378.


\bibitem{Ci}
{C.~Cibils:} {\em Review of D.~Happel: Hochschild cohomology of 
finite-dimensional algebras.  S\'eminaire d'alg\`ebre P. Dubreil et 
M.-P. Malliavin, Proc., Paris/Fr.  1987/88, Springer Lect.  Notes 
Math.  {\bf 1404} (1989), 108--126\/}; Math.  Reviews {\bf 91 b:16012} 
(1991).

\bibitem{Di}
{P.~A.~M.~Dirac:}
{\em The Fundamental Equations of Quantum Mechanics.\/}
Proc. Roy. Soc. A {\bf 109}, p. 642; reprinted in:
{\em Sources of Quantum Mechanics\/}, ed. by B.~L.~van der Warden, 
Dover Publications, Inc., New York 1968

\bibitem{Ga}
{P.~Gabriel}, {\em Auslander-Reiten sequences and representation-finite
algebras\/}, Lecture Notes in Mathematics, {\bf 831}, (Springer-Verlag, Berlin,
1980), 1--71.

\bibitem{Gas}
{S.~Gastamina, J.~A.~de la Pe\~{n}a, M.~I~ Platzeck, M.~J.~Redondo and
 S.~Trepode}, {\em Finite dimensional algebras
with vanishing Hochschild cohomology\/}, preprint.

\bibitem{Hap}
D.~Happel:
{\em
Hochschild cohomology of finite-dimensional algebras.\/}
S\'eminaire d'alg\`ebre P. Dubreil et M.-P. Malliavin, Proc.,
Paris/Fr. 1987/88, Lecture Notes in Mathematics, {\bf 1404}
(Springer, Berlin, 1989), 108--126.

\bibitem{Ig}
K.~Igusa:
{\em Notes on the no loop conjecture.\/}
J. Pure and Applied Algebra, {\bf 69} (1990), 161--176.
\bibitem{Len}
{H.~Lenzing:}
{\em Nilpotente Elemente in Ringen von endlicher globaler Dimension\/},
Math. Z., {\bf 108} (1969), 313--324.


\bibitem{Lod}
{J.~L.~Loday:}
{\em Cyclic Homology.\/}
Grundlehren der math. Wissenschaften, {\bf 301},
Springer-Verlag, New York 1992.

\end{thebibliography}
\end{document}